\documentclass[reqno,oneside,12pt]{amsart}
\usepackage{amssymb}
\begin{document}
\newtheorem{prop}{Proposition}
\newtheorem*{teo}{Theorem}
\title{Harmonic Oscillator in Characteristic $p$}
\author{Anatoly N. Kochubei}
\address{Institute of Mathematics,
Ukrainian National Academy of Sciences,
Tereshchenkivska 3, Kiev, 252601 Ukraine}
\email{ank@ank.kiev.ua}
\thanks{Partially supported by
the Ukrainian Fund of Fundamental Research (Grant 1.4/62)}
\date{}
\begin{abstract}
We construct an irreducible representation of the canonical
commutation relations by operators on a certain Banach space over
a local field of characteristic $p$. The Carlitz polynomials
forming the basis of the space are shown to be the counterparts
of the Hermite functions for this situation. The analogues of
coherent states are related to the Carlitz exponential.
\end{abstract}
\maketitle
\section{Introduction}

It was shown in [Ko] that the canonical and deformed commutation
relations admit representations by bounded operators in Banach
spaces over the field of $p$-adic numbers, and that various
objects related to those representations (analogues of the
Hermite functions, coherent states etc.)  coincide with
well-known special functions of $p$-adic analysis. The constructions
lead to the introduction of some operators (analogues of the
number operator) which possess eigenbases orthonormal in the
non-Archimidean sense. Note that no general concept of a ``Laplacian''
is known in $p$-adic analysis, and $p$-adic spectral theory (see
[B, Vi]) provides no tool to establish the ``Hermitian'' property of an
operator, other than to construct its eigenbasis. Thus the simple
difference operators introduced in [Ko] and related to the
additive and multiplicative structures of the $p$-adics, may be
viewed as Laplacians on model $p$-adic domains.

In this Letter we pursue the same line for another model example,
the case of a local field of characteristic $p$, that is of the
field $K$ of formal Laurent series with coefficients from the
Galois field $\mathbf F_q$. Here $p$ is a prime number, $q=p^\gamma
$, $\gamma \in \mathbf Z_+$. The construction is very simple again,
though quite different from the ones in [Ko], and involves
features typical for analysis over $K$ -- the basic space is
generated by $\mathbf F_q$-linear polynomials, the operators
themselves are only $\mathbf F_q$-linear. This time the well-known
Carlitz polynomials (see [C1, C2, G1])  appear as analogues of
the Hermite functions, and the ``coherent states'' (eigenfunctions
of the annihilation operator) are expressed via the Carlitz
exponential, the basic object in the arithmetic of function
fields [AT, G2, G3].

\section{Carlitz Basis}

Let us denote by $|\cdot |$  the non-Archimedean absolute value in
$K$; if $z\in K$,
$$
z=\sum \limits _{i=n}^\infty a_ix^i,\quad n\in \mathbf Z ,\ a_i\in
\mathbf F_q,
$$
then $|z|=q^{-n}$. Let $O=\{ z\in K: |z|\le 1\} $ be the ring of
integers in $K$; the ring $\mathbf F_q[x]$ of polynomials (in the
indeterminate $x$) with coefficients from $\mathbf F_q$ is dense in $O$.
Let $K_{ac}$ be the algebraic closure of $K$, with the absolute value
induced from $K$. We shall denote by $\overline{K}$ the
completion of $K_{ac}$. It is known [K\"u, R] that $\overline{K}$
is an algebraically closed field.

A function $\varphi :\ O\to \overline{K}$ is called
$\mathbf F_q$-linear if $\varphi (t_1+t_2)=\varphi (t_1)+\varphi
(t_2)$ and $\varphi (\beta t)=\beta \varphi (t)$ for any
$t,t_1,t_2\in O$, $\beta \in \mathbf F_q$. Let $X$ be a Banach
space (over $\overline{K}$) of all continuous $\mathbf F_q$-linear
$\overline{K}$-valued functions, with the supremum norm.

The Carlitz polynomials $f_i(t),\ i=0,1,2,\ldots $, are defined
as follows. Let $e_0(t)=t$,
\begin{equation}
e_i(t)=\prod \limits _{{m\in \mathbf F_q[x]\atop \deg m<i}}(t-m),\quad
i\ge 1.
\end{equation}
It is known [C1, G1] that
\begin{equation}
e_i(t)=\sum \limits _{j=0}^i(-1)^{i-j}\left[ i\atop j\right]
t^{q^j}
\end{equation}
where
$$
\left[ i\atop j\right]=\frac{D_i}{D_jL_{i-j}^{q^j}},
$$
the elements $D_i,\ L_i\in K$ are defined as
$$
D_i=[i][i-1]^q\ldots [1]^{q^{i-1}};\ L_i=[i][i-1]\ldots [1]\
(i\ge 1);\ D_0=L_0=1,
$$
and $[i]=x^{q^i}-x\in O$. Finally,
\begin{equation}
f_i(t)=D_i^{-1}e_i(t),\quad i=0,1,2,\ldots  .
\end{equation}
Since char $K=p$, we know that
$(t_1+t_2)^{q^j}=T_1^{q^j}+t_2^{q^j}$ for any $t_1,t_2\in K$; if
$\beta \in \mathbf F_q$ then $\beta ^{q^j}=\beta $. Hence, $f_i\in
X$, $i=0,1,\ldots $.

It was shown by Wagner [W] (see also [G1]) that $\{ f_i\} $ is a
basis  in $X$, that is any function $\varphi \in X$ admits a
unique representation as a uniformly convergent series
$$
\varphi =\sum \limits _{i=0}^\infty c_if_i, \quad c_i\in
\overline{K},\ c_i\to 0,
$$
with  $|c_i|\le 1$ for all $i$, if $\Vert \varphi \Vert \le 1$.

\begin{prop}
The basis $\{ f_i\} $ is orthonormal (in the non-Archimedean
sense), that is
$$
\Vert \varphi \Vert =\sup \limits _{i\ge 0}|c_i|, \quad \mbox{for
any } \varphi \in X.
$$
\end{prop}

The proof will follow immediately from the next, more general
proposition, which is of some independent interest.

According to [C2, W, G1], the basis $\{ f_i\} $ can be
complemented up to a basis of the space $C(O,\overline{K})$ of
all continuous $\overline{K}$-valued functions on $O$, as
follows. Let us write any natural number $j$ as
\begin{equation}
j=\sum \limits _{i=0}^\nu \alpha _iq^i,\quad 0\le  \alpha _i<q,
\end{equation}
and set
\begin{equation}
h_j(t)=\frac{G_j(t)}{\Gamma _j},\quad G_j(t)=\prod \limits
_{i=0}^\nu (e_i(t))^{\alpha _i},\quad  \Gamma _j=\prod \limits
_{i=0}^\nu D_i^{\alpha _i},
\end{equation}
so that $h_{q^i}=f_i,\ \ h_j$ is a polynomial of degree $j$.

\begin{prop}
The basis $\{ h_j\} $ is orthonormal.
\end{prop}

{\it Proof.} Let us start from another basis $\{ Q_j\} $ in
$C(O,\overline{K})$ whose orthonormality follows from the results
of [W,A]. Writing any natural number $j$ in the form (4), set
\begin{equation}
Q_0(t)\equiv 1,\quad Q_j(t)=\frac{P_j(t)}{P_j(m_j)},\ j\ge 1,
\end{equation}
where
$$
P_j(t)=(t-m_0)(t-m_1)\ldots (t-m_{j-1}),\ \ m_j=a_{\alpha
_0}+a_{\alpha _1}x+\cdots +a_{\alpha _\nu }x^\nu ,
$$
and $a_k$ are the elements of $\mathbf F_q=\{ a_0,\ldots ,a_{q-
1}\}$, $a_0=0$, $a_1=1$ (we do not show the dependence of $\nu,
\alpha _1,\ldots ,\alpha _\nu $ on $j$).

Polynomials $\{ Q_j\} _{j\le n}$  form a basis in the space of
all polynomials with degrees $\le n$. In particular, one can
write
$$
h_n(t)=\sum \limits _{i=0}^nc_{ni}Q_i(t),\quad n=0,1,2,\ldots .
$$
Since $h_n(t)\in \mathbf F_q[x]$ for all $t\in \mathbf F_q[x]$ (see
[C2]), we have $\Vert h_n\Vert \le 1$, whence $|c_{ni}|\le 1$ for
all $n,i$.

Now, by a general result from [Ve], in order to prove the
proposition, it suffices to show that
$$
|c_{nn}|=1, \quad n=0,1,2,\ldots .
$$

It is seen from (1), (3), (5) and (6) that the leading
coefficient of $h_n(t)$ equals $\Gamma _n^{-1}$ while the leading
coefficient of $Q_n(t)$ is $(P_n(m_n))^{-1}$ so that
$c_{nn}=\Gamma _n^{-1}P_n(m_n)$.

By definition,
$$
D_n=x^{1+q+\cdots +q^{n-1}}\left( x^{q^n-1}-1\right) \left(
x^{q^{n-1}-1}-1\right) ^q\ldots \left( x^{q-1}-1\right) ^{q^{n-
1}},
$$
so that $\Gamma _n=x^{l_n}S_n(x)$ where $S_n\in \mathbf F_q[x],\
S_n(0)\ne 0$,
$$
l_n=\alpha _1+(1+q)\alpha _2+\cdots +(1+q+\cdots +q^s)\alpha _s,
$$
and $\alpha _1,\ldots \alpha _s$ are taken from the $q$-adic
expansion
\begin{equation}
n=\alpha _0+\alpha _1q+\cdots +\alpha _sq^s,\quad 0\le  \alpha
_j<q.
\end{equation}
We have $|\Gamma _n|=q^{-l_n}$.

On the other hand,
$$
P_n(m_n)=(m_n-m_0)(m_n-m_1)\ldots (m_n-m_{n-1})=x^{\kappa
_n}T_n(x),
$$
where $T_n\in \mathbf F_q[x],\ T_n(0)\ne 0$, while the number
$\kappa _n$ can be computed as follows. Consider a difference
$m_n-m_{n'}$ with $n'<n$, and write
$$
n'=\alpha '_0+\alpha '_1q+\cdots +\alpha '_sq^s,\ \
m_{n'}=a_{\alpha '_0}+a_{\alpha '_1}x+\cdots +a_{\alpha '_s }x^s.
$$
The quantity of the differences $m_n-m_{n'}$ divisible by $x$
equals the quantity of natural numbers of the form
$\alpha _0+\alpha '_1q+\cdots +\alpha '_sq^s$ which are less or
equal  $n$ (here $\alpha _0$ is just the same as in (7)), that is
the quantity of $s$-tuples  $(\alpha '_1,\ldots \alpha '_s)$ such
that
$$
\alpha '_1+\alpha '_2q+\cdots +\alpha '_sq^{s-1}\le
\alpha _1+\alpha _2q+\cdots +\alpha _sq^{s-1}.
$$
Of course, such a quantity equals
$\alpha _1+\alpha _2q+\cdots +\alpha _sq^{s-1}$.

Similarly, the quantity of the differences divisible by $x^2$ is
$\alpha _2+\alpha _3q+\cdots +\alpha _sq^{s-2}$, so that the
number of differences divisible by $x$ and not divisible by $x^2$
equals to
$$
(\alpha _1+\alpha _2q+\cdots +\alpha _sq^{s-1})-
(\alpha _2+\alpha _3q+\cdots +\alpha _sq^{s-2}).
$$

Continuing this reasoning we find that
$$
\begin{array}{l}
\kappa _n=(\alpha _1+\alpha _2q+\cdots +\alpha _sq^{s-1})-
(\alpha _2+\alpha _3q+\cdots +\alpha _sq^{s-2})\\ \qquad {}+
2(\alpha _2+\alpha _3q+\cdots +\alpha _sq^{s-2})-
(\alpha _3+\alpha _4q+\cdots +\alpha _sq^{s-3})+\cdots \\ \qquad +(s-
1)((\alpha _{s-1}+\alpha _sq)-\alpha _s)+s\alpha _s\\ \qquad {}=
(\alpha _1+\cdots +\alpha _sq^{s-1})+(\alpha _2+\cdots +\alpha _sq^{s-2})
+\cdots +(\alpha _{s-1}+\alpha _sq)\\ \qquad +\alpha _s=l_n
\end{array}
$$
which means that $|c_{nn}|=1.\quad \square $

\section{Carlitz Exponential}

The Carlitz exponential is defined by the power series
\begin{equation}
e_C(z)=\sum \limits _{j=0}^\infty \frac{z^{q^j}}{D_j}.
\end{equation}
Since
$$
|D_j|=q^{-1}(q^{-1})^q\ldots (q^{-1})^{q^{j-1}}=q^{-\frac{q^j-
1}{q-1}}=\left( q^{-1/(q-1)}\right) ^{q^j-1},
$$
the series (8) is convergent if $z\in \overline{K},\ |z|<q^{-1/(q-1)}$.

We shall need some relations involving $e_C$, which coincide
formally with the ones well known in number theory (see e.g. [G1]
or [AT]), though in the literature the series (8) is usually
considered not over $\overline{K}$, but rather over the $\infty
$-adic completion of the field of rational functions (with the
field of constants $\mathbf F_q)$.

Consider a function
\begin{equation}
\rho (\zeta )=\sum \limits _{n=0}^\infty (-1)^n\frac{\zeta
^{q^n}}{L_n}.
\end{equation}
It is known [AT] that the functions $e_C$ and $\rho $ are inverse
to each other as formal power series. Since $|L_n|=q^{-n}$, the
series (9) is convergent if $|\zeta |<1$, and
$$
\rho (\zeta )|\le \sup \limits _{n\ge 0}q^n|\zeta |^{q^n}.
$$
The function $\psi _\zeta (s)=s|\zeta |^s$ is decreasing for $s>-
(\log |\zeta |)^{-1}$; if $|\zeta |<q^{-1/(q-1)}$ then $\psi
_\zeta $ decreases for $s>(q-1)(\log q)^{-1}$. In particular,
$\psi _\zeta (q^n)\le \psi _\zeta (q)$, $n\ge 1$. Hence, $\rho
(\zeta )|\le \max (|\zeta |,q|\zeta |^q)$, and we find that
$$
|\rho (\zeta )|<q^{-1/(q-1)}\quad \mbox{if }\ |\zeta |<q^{-1/(q-
1)},
$$
which implies the identity
\begin{equation}
e_C(\rho (\zeta ))=\zeta ,\quad |\zeta |<q^{-1/(q-1)}.
\end{equation}

For any fixed $z\in \overline{K}$ such that $|\zeta |<q^{-1/(q-
1)}$, consider the function $w_z\in X$ of the form
$$
w_z(t)=e_C(tz),\quad t\in O.
$$
Let us find its expansion with respect to the Carlitz basis.
Below we shall use the difference operators
$$
\Delta \varphi (t)\equiv \Delta ^{(1)}\varphi (t)=\varphi (xt)-
x\varphi (t);
$$
$$
\Delta ^{(i)}\varphi (t)=\Delta ^{(i-1)}\varphi (xt)-
x^{q^{i-1}}\Delta ^{(i-1)}\varphi (t), \quad i\ge 2.
$$

\begin{prop}
The function $w_z$ can be expanded as
\begin{equation}
w_z(t)=\sum \limits _{n=0}^\infty (e_C(z))^{q^n}f_n(t).
\end{equation}
\end{prop}

{\it Proof.} Let
$$
e_C^{(N)}(z)=\sum \limits _{j=0}^N \frac{z^{q^j}}{D_j};\quad
w_{N,z}(t)=e_C^{(N)}(tz).
$$
It is clear that $w_{N,z}(t)\longrightarrow w_z(t)$ uniformly
with respect to $t\in O$. On the other hand (see [G1]),
$$
w_{N,z}(t)=\sum \limits _{n=0}^Nb_n^{(N)}(z)f_n(t),
$$
where $b_n^{(N)}(z)=\Delta ^{(n)}e_C^{(N)}(z),\ n\le N,\ \Delta ^{(0)}=I$
(the identity operator).

Using (8), the relation $D_i=[i]D_{i-1}^q$, and employing repeatedly
the binomial  relation for a field of characteristic $p$, we find by
induction that
$$
b_n^{(N)}(z)=\left( e_C^{(N-n)}(z)\right) ^{q^n},\quad n\ge 0.
$$

It follows from orthonormality of the Carlitz basis that for $N\to \infty $
the coefficient $b_n^{(N)}(z)$ converges to the $n$-th coefficient of the
expansion of the function $w_z$. Thus we come to (11). $\quad \square $

\section{Creation and Annihilation Operators}

The operator $A:\ X\to X$ is called $\mathbf F_q$-linear if $A(u+v)=
Au+Av$ and $A(\beta u)=\beta Au$ for any  $u,v\in X$, $\beta \in
\mathbf F_q$. The simplest example is the operator $R_qu=u^q$. Its
inverse $\sqrt[q]{\ }$ defined on the set $\{ u^q\ :\ u\in X\} $
possesses similar properties (recall that $q=p^\gamma $, so that the
$q$-th root, if it exists, is unique in a field of characteristic
$p$).

Now we can introduce our oscillator-like model. The formula (13) below
can create an impression that the notation $[i]$ was invented
deliberately in order to emphasize the analogy with the conventional
quantum mechanics. In reality this notation was proposed by Carlitz
in 1935!

Let
$$
a^+=R_q-I,\quad a^-=\sqrt[q]{\ }\circ \Delta .
$$

\begin{teo}

{\rm (i)} $a^+$ and $a^-$ are continuous $\mathbf F_q$-linear operators
on $X$,
\begin{equation}
a^-a^+-a^+a^-=[1]^{1/q}I.
\end{equation}
{\rm (ii)} The operator $a^+a^-$ possesses the othonormal eigenbasis
$\{ f_i\} $,
\begin{equation}
(a^+a^-)f_i=[i]f_i,\quad i=0,1,2,\ldots ;
\end{equation}
$a^+$ and $a^-$ act upon the basis as follows:
\begin{equation}
a^+f_{i-1}=[i]f_i,\ \ a^-f_i=f_{i-1},\ i\ge 1;\ a^-f_0=0.
\end{equation}
{\rm (iii)} The equation
\begin{equation}
a^-u=\lambda u
\end{equation}
has solutions (``coherent states'') for any $\lambda \in
\overline{K}$; if $\lambda \ne 0$, each solution can be written as
\begin{equation}
u(t)=\lambda ^{-q/(q-1)}\sum \limits _{n=0}^\infty
c^{q^n}f_n(t),\quad c\in \overline{K},\ |c|<1,
\end{equation}
for some value of the $(q-1)$-th root, and conversely, $a^-
u=\lambda u$ for the function (16). If in (16) $|c|<q^{-1/(q-1)}$
then
\begin{equation}
u(t)=\lambda ^{-q/(q-1)}e_C(tz),\quad z=\rho (c).
\end{equation}
In particular, if $q\ne 2$ then every function (16) with $c\in K$
takes the form (17).
\end{teo}

{\it Proof. }
(i) The operator $\Delta $ is linear and transforms any function
$\varphi \in X$ into the $q$-th power of some function  from $X$.
Indeed,
$$
\Delta e_i=\frac{D_i}{D_{i-1}^q}e_{i-1}^q,\ i\ge 1;\ \ \Delta e_0=0
$$
(see [G1]) whence
\begin{equation}
\Delta f_i=\left\{
\begin{array}{rl}
f_{i-1}^q, &\ \ i\ge 1;\\
0, &\ \ i=0.
\end{array}\right.
\end{equation}

If
$$
\varphi (t)=\sum \limits _{i=0}^\infty c_if_i(t),\quad c_i\to 0,
$$
then (recall that char $\overline{K}=p$)
$$
\Delta \varphi (t)=\sum \limits _{i=1}^\infty c_if_{i-
1}^q(t)=\left( \sum \limits _{i=1}^\infty c_i^{1/q}f_{i-
1}(t)\right) ^q.
$$
Thus $a^-$ is correctly defined, $\mathbf F_q$-linear, and
\begin{equation}
a^-\varphi (t)=\sum \limits _{i=1}^\infty c_i^{1/q}f_{i-
1}(t),
\end{equation}
so that $\Vert a^-\varphi \Vert =\sup \limits _{i\ge
1}|c_i^{1/q}|\le \Vert \varphi \Vert $ which implies continuity
of $a^-$. Similar properties of $a^+$ are obvious.

Simple calculation yeilds the formulas
$$
(a^+a^-\varphi )(t)=\Delta \varphi (t)-(\Delta \varphi
(t))^{1/q};
$$
$$
(a^-a^+\varphi )(t)=(\Delta \varphi ^q(t))^{1/q}-(\Delta \varphi
(t))^{1/q}.
$$
Subtracting we get
$$
(a^-a^+-a^+a^-)\varphi (t)=(x-x^{1/q})\varphi (t)=(x^q-
x)^{1/q}\varphi (t),
$$
and we come to (12).

(ii) The formula (13) is a consequence of (14); the latter
follows from (18) and the identities
$$
e_i=e_{i-1}^q-D_{i-1}^{q-1}e_{i-1},\quad D_i=[i]D_{i-1}^q
$$
(see [G1]).

(iii)  Let $a^-u=\lambda u$,
$$
u(t)=\sum \limits _{n=0}^\infty c_nf_n(t),\quad c_n\to 0.
$$
It follows from (19) and the uniqueness of the expansion that
$$
c_{n+1}^{1/q}=\lambda c_n,\quad n=0,1,\ldots ,
$$
whence
$$
c_n=\lambda ^{q^n+q^{n-1}+\cdots +q}c_0^{q^n}=\mu ^{-1}(c_0\mu
)^{q^n},\quad n=1,2,\ldots ,
$$
where $\mu =\lambda ^{q/(q-1)}$. Since $c_n\to 0$, we have
$|c_0\mu |<1$, and we obtain the representation (16) with
$c=c_0\mu $.

The converse statement follows easily from the identity (18).

The representation (17) is a direct consequence of (16) and the
results of Sect. 3. If $q\ne 2$, $c\in K$, $|c|<1$, then $|c|\le q^{-
1}<q^{-1/(q-1)}$, so that in this case (16) is equivalent to
(17). $\quad \square $

\end{document}